# Arithmetical Properties of Laplacians of Graphs

Dino Lorenzini

January 6, 1999

Let $M \in M_n(\mathbb{Z})$ denote any matrix. Thinking of $M$ as a linear map $M : \mathbb{Z}^n \to \mathbb{Z}^n$, we denote by $\text{Im}(M)$ the $\mathbb{Z}$-span of the column vectors of $M$. Let $e_1, \ldots, e_n$, denote the standard basis of $\mathbb{Z}^n$, and let $E_{ij} := e_i - e_j$, $(i \neq j)$. In this article, we are interested in the group $\mathbb{Z}^n/\text{Im}(M)$, and in particular in the elements of this group defined by the images of the vectors $E_{ij}$ under the quotient $\mathbb{Z}^n \to \mathbb{Z}^n/\text{Im}(M)$. Most of this article is devoted to the study of the case where $M$ is the laplacian of a graph.

Let $G$ be any connected graph with $n$ vertices and $m$ edges. (In this paper, all graphs are assumed to be connected unless stated otherwise, and any edge of $G$ has two distinct endpoints.) Choose an ordering of the vertices of $G$. If $i \neq j$, then let $c_{ij}$ denote the number of edges linking the vertex $v_i$ with the vertex $v_j$. Let us say that a graph is *simple* if $c_{ij} \leq 1$, for all $i \neq j$. Let $-c_{ii} := \sum_{j \neq i} c_{ij}$. The integer $c_{ii}$ is called the degree of $v_i$, and may also be denoted by $d_i$. Let $M(G) := ((c_{ij}))_{1 \leq i,j \leq n}$. The matrix $-M(G)$ is called the laplacian of $G$. For practical reasons, we prefer to work with $M(G)$ rather than with $-M(G)$ since the former has fewer negative coefficients.

It is well-known that the group $\mathbb{Z}^n/\text{Im}(M(G))$ is the product of $\mathbb{Z}$ with a finite abelian group that we shall denote by $\Phi(G)$. In fact, the group $\mathbb{Z}^n/\text{Im}(M(G))$ can be described using any diagonal matrix $\text{diag}(m_1, \ldots, m_{n-1}, 0)$ which is row and column equivalent over $\mathbb{Z}$ to $M$:

$$\mathbb{Z}^n/\text{Im}(M(G)) \cong (\prod_{i=1}^{n-1} \mathbb{Z}/m_i\mathbb{Z}) \times \mathbb{Z}.$$

Thus the study of the group $\mathbb{Z}^n/\text{Im}(M(G))$ is closely related to the study of the Smith Normal Form of $M(G)$. The kernel of $M(G)$ is generated by the transpose $V$ of the vector $(1, \ldots, 1)$. Since the vectors $E_{ij}$ form a basis of the orthogonal space to the vector $V$, we find that the set $\{E_{ij}, i \neq j\}$ is a set of generators for $\Phi(G)$.

The relationships between the group $\mathbb{Z}^n/\operatorname{Im}(M(G))$ and the geometry of the graph $G$ are not yet well understood. It is known that the order of $\Phi(G)$ is equal to the number $\kappa(G)$ of spanning trees of $G$. Let $\beta(G) := m - n + 1$ denote the number of independent cycles of $G$. It is shown in [Lor1], p. 281, that $\Phi(G)$ can be generated by $\beta(G)$ elements. The number of generators is also bounded by $n - 1 - \operatorname{diam}(G)$ in [GMW], Cor. 1. For an open question regarding a possible relationship between the geometry of $G$ and the group $\Phi(G)$, see [Lor2].

We shall say that *a pair of vertices* $\{v_i, v_j\}$ *on a graph* $G$ *has order* $h > 0$ if there exists a vector $S \in \mathbb{Z}^n$ with $^tS = (s_1, \ldots, s_n)$ such that $MS = hE_{ij}$ and $\gcd(s_1 - s_n, \ldots, s_{n-1} - s_n) = 1$. Thus a pair of vertices $\{v_i, v_j\}$ of $G$ has order $h$ if and only if the image of $E_{ij}$ has order $h$ in $\mathbb{Z}^n/\operatorname{Im}(M(G))$. Note that if $S$ and $S'$ are such that $MS = MS'$, then $^tS' = {}^tS + k(1, \ldots, 1)$ for some integer $k$.

The smallest graph $G_h$ with a pair of vertices of order $h$ is such that

$$M(G_h) = \begin{pmatrix} -h & h \\ h & -h \end{pmatrix}.$$

It is an interesting open question to describe a smallest *simple* graph $G$ with a pair of order $h$. We describe in the next section several operations on graphs which, given a graph with a pair of order $h$, produce a new graph with a pair of order $h$. Our main Theorem 2.1 states that every graph $G$ with a pair of order $h$ can be obtained by starting with the graph $G_h$ and performing a sequence of such operations.

It would be very interesting to have an explicit graph-theoretical algorithm for computing the order of a given pair of vertices $\{v_i, v_j\}$ on a graph $G$. In this paper, we provide an explicit criterion for determining whether a given pair $\{v_i, v_j\}$ has order 1 or 2 (see 2.3, 2.4). Let us state in this introduction the following related result linking the geometry of the graph and properties of the group $\mathbb{Z}^n/\operatorname{Im}(M(G))$.

Let us call a graph *multiply connected* if, given any edge $e$ of $G$, the graph $G \setminus \{e\}$ is connected. Let us call a matrix $M \in M_n(\mathbb{Z})$ *spread* if the quotient map $\mathbb{Z}^n \to \mathbb{Z}^n/\operatorname{Im}(M)$ is injective when restricted to the set $\{e_1, \ldots, e_n\}$. In other words, $M$ is spread if none of the images of the vectors $E_{ij}$ in $\mathbb{Z}^n/\operatorname{Im}(M)$ has order 1.

**Proposition 2.3** *Let $G$ be a connected graph. Then $G$ is multiply connected if and only if the laplacian of $G$ is spread.*

In the last section of this article, we investigate the properties of the groups $\mathbb{Z}^n/\operatorname{Im}(M - \mu\operatorname{Id})$, where $\mu$ is any integer. In the case of a laplacian $M$, we show in Theorem 3.1 that the matrix $M - \mu\operatorname{Id}$ is spread if $\mu \notin [0, n+1]$. A similar statement holds for any semidefinite positive quadratic form (3.4).

Our initial motivation for studying groups of the type $\mathbb{Z}^n/\operatorname{Im}(M)$ is found in arithmetic geometry. Let $K$ be any field. Let $X/K$ be any smooth proper geometrically connected curve. Associated to this curve in a functorial way is a group variety $A/K$, called the jacobian of $X$. When $X$ has a $K$-rational point, say $P_0$, then one can define a map of algebraic varieties $\varphi : X \to A$. Call $\varphi_K : X(K) \to A(K)$ the restriction of this map to the sets of $K$-rational points. Then $\varphi_K$ is simply $P \mapsto P - P_0$. Since $A(K)$ is a group, it is natural to wonder whether $P - P_0$ generates a finite group in $A(K)$. If it does, then we would like to know the order of this group. These questions are usually very hard to answer.

When $K$ is a field with a discrete valuation $v$, with ring of integers $\mathcal{O}_K$, and with maximal ideal $(\pi)$, it is possible to reduce modulo $(\pi)$ both the curve $X$ and its jacobian $A$. One obtains using standard techniques a group homomorphism, red : $A(K) \to \Phi$, where the group $\Phi$ is a finite abelian group called the group of components of the Néron model of $A/K$. Thus, information on $P - P_0$ can be obtained by understanding first red$(P - P_0)$. Associated to the reduction of the curve $X$ is an intersection matrix $M$ and an arithmetical graph $(G, M, R)$ as in [Lor3]. The group $\Phi$ is isomorphic to the torsion subgroup of $\mathbb{Z}^n/\operatorname{Im}(M)$. The image of a point $P - P_0$ in $A(K)$ under the reduction map corresponds to the class in $\mathbb{Z}^n/\operatorname{Im}(M)$ of a vector of the form $E_{ij}$.

When the curve $X$ has semistable reduction, the associated intersection matrix $M$ is simply, up to a sign, the laplacian of a connected graph $G$. Thus the results of this article immediately apply to this situation. The general case where $M$ is any intersection matrix is more delicate. Further applications of the results of this article to this case can be found in the forthcoming article [Lor4]. No further reference to arithmetic geometry is made in this article, except for an application to modular curves in 2.5.

The author would like to thank H.-G. Evertz, A. Granville, and C. Pomerance, for helpful discussions.

# 1 Constructions of $h$-marked graphs

Let us call *$h$-marked graph* a graph $G$ with a pair of vertices $\{v_i, v_j\}$ of order $h$. The data of a vector $S$ such that $MS = hE_{ij}$ for some $i \neq j$ (and $\gcd(s_1 - s_n, \ldots, s_{n-1} - s_n) = 1$) will be called a *marking* on $G$. The integer $s_k$ in ${}^tS = (s_1, \ldots, s_n)$ will be called the *weight* of the vertex $v_k$. Given a marking $S$ on a graph $G$, ${}^tS' := {}^tS + \alpha(1, \ldots, 1)$ also produces a marking of the vertices of $G$, and we shall always consider two such markings as equivalent.

**Lemma 1.1** *Let $\{v_i, v_j\}$ be a pair of vertices of order $h$ in $G$, with associated vector ${}^tS = (s_1, \ldots, s_n)$. Then $s_i = \min\{s_\ell, \ell = 1, \ldots, n\}$, and $s_j = \max\{s_\ell, \ell = 1, \ldots, n\}$.*

*Proof:* Let $s := \min\{s_\ell, \ell = 1, \ldots, n\}$. Since ${}^tS$ is not a multiple of $(1, \ldots, 1)$, we can find a vertex $v_k$ such that $s_k = s$ and there exists a vertex $v_p$ linked to $v_k$ with $s_p > s_k$. Then

$$\sum_{\ell=1}^n c_{k\ell} s_\ell = \sum_{\ell=1}^n c_{k\ell}(s_\ell - s_k) \geq \sum_{c_{k\ell} \neq 0} (s_\ell - s_k) > 0.$$

Since the vector $MS$ has a positive entry only in its $i$-th row, we find that $k = i$. A similar argument shows that $s_j = \max\{s_\ell, \ell = 1, \ldots, n\}$.

**Example 1.2** A cycle of $n$ vertices can be $n$-marked as follows:

where ${}^tS = (0, 1, \ldots, n-1)$. The reader will verify that the pair $\{v_1, v_n\}$ has order $n$.

We shall usually not number explicitly the vertices of a graph as on the picture on the left above, but use a picture as on the right above to explicitly indicate an $h$-marking of the graph. When drawing an $h$-marked graph, we will, when no confusion may occur, call the vertex of lowest weight $v_1$, and the vertex of highest weight $v_n$. The pair $\{v_1, v_n\}$ has order $h$, with $h$ determined as follows:

$$h = \sum_{\text{all } v \text{ adjacent to } v_1} (\text{weight of } v - \text{weight of } v_1)(\# \text{ edges between } v \text{ and } v_1).$$

**Example 1.3** We list below a few examples of small graphs with low order $h$. In each case, take $v_1$ to be the vertex of weight 0, and $v_n$ to be the vertex of highest weight. Then $\{v_1, v_n\}$ is a pair of order $h$.

Let us now describe some constructions/operations on graphs which, when applied to a given $h$-marked graph, produce a new $h$-marked graph. We shall show in Theorem 2.1 that any $h$-marked graph $G$ can be obtained from the graph $G_h$ using a sequence of operations described below.

**Construction 1.4** *Adding a graph to a given $h$-marked graph $G$.*

Pick any vertex on the given graph $G$ and number it $v_n$. The other vertices of $G$ are denoted by $v_1, \ldots, v_{n-1}$. Let $G'$ be any other graph. Pick any vertex of $G'$ and call it $w_1$. The other vertices of $G'$ are then denoted by $w_2, \ldots, w_{n'}$. Let $G''$ denote the union of the graphs $G$ and $G'$, with $v_n$ identified to the vertex $w_1$. Let $\{v_i, v_j\}$ denote the given pair of order $h$ in $G$, with associated vector ${}^tS = (s_1, \ldots, s_n)$. We shall number the vertices of $G''$ by $v_1, \ldots, v_n = w_1, w_2, \ldots, w_{n'}$. Then $G''$ is a new $h$-marked graph as follows: associate to the pair $\{v_i, v_j\}$ the vector

$$ {}^tS'' := (s_1, \ldots, s_n, s_n, \ldots, s_n). $$

In other words, the weight of $w_i$ in $G''$ is $s_n$ and the weight of $v_k$ is $G''$ is $s_k$. We shall say that the $h$-marked graph $G''$ is obtained by *adding the graph $G'$ to the $h$-marked graph $G$ at the vertex $v_n$*.

**Construction 1.5** *Gluing vertices of same weight on a given $h$-marked graph $G$.*

Let $\{v_1, v_n\}$ be the given pair of order $h$ on $G$. Let $v_i$ and $v_j$ be two vertices different from $v_1$ and $v_n$. Let ${}^tS = (s_1, \ldots, s_n)$ denote the vector associated to the pair $\{v_1, v_n\}$. Assume that $s_i = s_j$, and that $c_{ij} = 0$. Consider the graph $G'$ obtained by identifying the vertices $v_i$ and $v_j$ in a single vertex $v$. The vertices of $G'$ are then

$$ (\{v_1, \ldots, v_n\} \setminus \{v_i, v_j\}) \sqcup \{v\}. $$

The pair $\{v_1, v_n\}$ is of order $h$ in $G'$ when $v$ is given the weight $s_i$, and $v_k$ the weight $s_k$, for all $k \neq i, j$.

We shall say that the $h$-marked graph $G'$ is obtained *from $G$ by gluing $v_i$ to $v_j$*. In Example 1.3, f) is obtained from e) and h) is obtained from i) by gluing.

As mentioned in the introduction, the smallest $h$-marked graph $G_h$ has $\begin{pmatrix} -h & h \\ h & -h \end{pmatrix}$ as laplacian. To describe multiple edges using a diagram, it is convenient to introduce the following notation. The symbols in a graph $G$ will denote the fact that $v_i$ and $v_j$ are linked in $G$ by $e$ distinct edges. The graph $G_h$ is then represented as follows:

**Construction 1.6** *Thickening an edge of an $h$-marked graph $G$.*

Let $v$ and $v'$ be two vertices in $G$ linked by $e$ edges, with weights $s$ and $s'$. In other words, $G$ contains a diagram of the form

Assume that $s > s'$. Let $G'$ denote the graph obtained from $G$ as follows. Remove the $e$ edges between $v$ and $v'$ and introduce $s - s' - 1$ new vertices $w_1, \ldots, w_{s-s'-1}$, linked to $v$ and $v'$ as follows:

To turn $G'$ into an $h$-marked graph, set the weight of $w_i$ to be $s' + i$, and keep the weight of $v_k$ to be $s_k$. The original pair $\{v_i, v_j\}$ on $G$, which may contain one or both of the vertices $v$ and $v'$, is still a pair of order $h$ in the new graph $G'$. We shall say that $G'$ is obtained from $G$ *by thickening the edge(s) between $v$ and $v'$*. In Example 1.3, one builds c) from b), and g) from f), by thickening an edge.

To simplify the formula introduced below, let us assume that the marking $S$ on $G$ is such that the weight $s'$ of $v'$ is equal to 0. (If necessary, add $(-s', \ldots, -s')$ to ${}^t S$.) Recall that $|\Phi(G)|$ is the number of spanning trees of $G$. Let $G_1$ denote the graph obtained from $G$ by removing the $e$ edges linking

$v$ to $v'$. Let $G_2$ denote the graph obtained from $G_1$ by gluing $v$ to $v'$. It is well known that
$$|\Phi(G)| = |\Phi(G_1)| + e|\Phi(G_2)|.$$
Using induction and a similar type of argument for $G'$, we find that

$$\begin{aligned}
|\Phi(G')| &= |\Phi(G_1)|(es)^{s-1}s + |\Phi(G_2)|(es)^s \\
&= s(es)^{s-1}(|\Phi(G_1)| + e|\Phi(G_2)|) \\
&= s(es)^{s-1}|\Phi(G)|.
\end{aligned}$$

In view of this nice formula for $|\Phi(G')|$, it is natural to wonder whether $\Phi(G)$ is a direct summand of $\Phi(G')$. The following two lemmas provide some information on the structure of the group $\Phi(G')$.

**Lemma 1.7** *The group $\Phi(G')$ contains $(\mathbb{Z}/es\mathbb{Z})^{s-2}$ as a direct summand.*

*Proof:* The pair of vertices $v'$ and $v$ in $G$ gives us a laplacian $M$ and a relation $MS = hE_{ij}$ of the form:

$$\begin{pmatrix} & \vdots & \vdots & \\ \cdots & -r' & e & \cdots \\ \cdots & e & -r & \cdots \\ & \vdots & \vdots & \end{pmatrix} \begin{pmatrix} \vdots \\ s' = 0 \\ s \\ \vdots \end{pmatrix} = hE_{ij}.$$

After the thickening, the vertices $v', w_1, \ldots, w_{s-1}, v$, of $G'$ give us

$$\begin{pmatrix} & \vdots & & & & & \\ \cdots & -r' + e - es & es & 0 & & & \\ & es & -2es & es & & & \\ & 0 & es & \ddots & \ddots & & \\ & & & \ddots & -2es & es & \\ & & & & es & -r + e - es & \cdots \\ & & & & & \vdots & \end{pmatrix} \begin{pmatrix} \vdots \\ 0 \\ 1 \\ \vdots \\ s-1 \\ s \\ \vdots \end{pmatrix}.$$

Add all columns to the $v'$-column to get a zero column. Then add all rows to the $v'$-row to get a zero $v'$-row. The central $(s-1) \times (s-1)$ submatrix is of the form
$$\begin{pmatrix} -2es & es & & \\ es & -2es & es & \\ & es & \ddots & \ddots \\ & & \ddots & -2es \end{pmatrix}.$$

An obvious row and column reduction over $\mathbb{Z}$ of this matrix produces a matrix

$$\begin{pmatrix} 0 & \cdots & \cdots & 0 & * \\ es & 0 & & & 0 \\ 0 & es & \ddots & & \vdots \\ \vdots & \ddots & \ddots & \ddots & 0 \\ & & & es & 0 \end{pmatrix}.$$

After permutations of the rows, we get a matrix of the form

$$\begin{pmatrix} es & 0 & \cdots & \cdots & 0 \\ 0 & es & \ddots & & \vdots \\ \vdots & \ddots & \ddots & & \vdots \\ & & & es & 0 \\ 0 & \cdots & \cdots & 0 & * \end{pmatrix}.$$

Thus, the group $\mathbb{Z}^n/\mathrm{Im}(M(G'))$ contains a direct summand of the form $(\mathbb{Z}/es\mathbb{Z})^{s-2}$.

**Lemma 1.8** *Let $G$ be an $h$-marked graph, with $\{v_1, v_n\}$ a pair of order $h$. Let ${}^tS = (0, s_2, \ldots, s_n = s)$ be its associated vector. Assume that $v_1$ and $v_n$ are linked by exactly one edge. Let $w_1, \ldots, w_{s-1}$ denote the vertices introduced in the thickening $G'$ of $G$ at the edge between $v_1$ and $v_n$. Then the pair $\{v_1, w_1\}$ has order $hs^2$ in $G'$. More precisely, drawn below is an $hs^2$-marking on $G'$. In this marking, the weight of $v_i$ is $s_i s$, and the weight of $w_i$ is $si + h(s-i)$.*

*Proof:* Left to the reader. The graph $G$ contains

and $G'$ contains

**Construction 1.9** *Removing from a given $h$-marked graph $G$ the edges of a subgraph.*

Let $\{v_i, v_j\}$ be the given pair of order $h$ on $G$, and let ${}^t S = (s_1, \ldots, s_n)$ be the associated vector. Consider a subgraph $G'$ of $G$, all of whose vertices have the same fixed weight $s_k$. Remove from $G$ all edges of $G'$, and call the resulting graph $\tilde{G}$. If $\tilde{G}$ is connected, then it is also $h$-marked with $S$. Assume now that $\tilde{G}$ is not connected.

**Lemma 1.10** *Both vertices $v_i$ and $v_j$ belong to the same connected component of $\tilde{G}$.*

*Proof:* Without loss of generality, we may assume that removing a single edge between two vertices $v$ and $v'$ of same weight $s$ disconnects the given connected graph $G$. Thus, the laplacian of $G$ is of the form

$$\begin{pmatrix} & & * & \vdots & & & \\ & & * & 0 & & & \\ * & * & -r & 1 & 0 & \cdots & \\ \cdots & 0 & 1 & -r' & * & * & \\ & & 0 & * & & & \\ & & \vdots & * & & & \end{pmatrix} \begin{pmatrix} * \\ * \\ s \\ s \\ * \\ * \end{pmatrix} = h E_{ij}.$$

Removing the edge between $v$ and $v'$ produces two new graphs $G_1$ and $G_2$. The graph $G_1$ has for laplacian the slightly modified top left corner of the laplacian of $G$:

$$M_1 := \begin{pmatrix} & & * \\ & & * \\ * & * & -r+1 \end{pmatrix}.$$

Suppose that $v_i$ is a vertex of $G_1$, but that $v_j$ is not. Then we find that

$$M_1 \begin{pmatrix} * \\ * \\ s \end{pmatrix} = h\, e_i, \text{with } h \neq 0.$$

This is a contradiction, since the column of the matrix $M_1$ are perpendicular to $(1, \ldots, 1)$ and $h e_i$ is not. Thus, Lemma 1.10 is proved.

Returning to our construction, we let $G'''$ denote the unique connected component of $\tilde{G}$ that contains both $v_i$ and $v_j$. By construction, $G'''$ is an $h$-marked graph when a vertex of $G'''$ is given as weight its weight in $G$. We shall call $G'''$ *the graph obtained by removing the subgraph $G'$ from the graph $G$*. In Example 1.3, the graph e) is obtained from d) by removing an edge.

**Construction 1.11** *Adding an edge to a given h-marked graph $G$.*

Let $G$ be an $h$-marked graph with a pair $\{v_i, v_j\}$ and associated vector $^tS$. Let $v_k$ and $v_\ell$ be two vertices with $s_k = s_\ell$. Let $G'$ be the graph obtained from $G$ by adding $c$ edges between $v_k$ and $v_\ell$. Then the pair $\{v_i, v_j\}$ is also of order $h$ on $G'$. The reader may check this fact directly, or may use Construction 1.4 to add a graph $G_c$ to $G$ at the vertex $v_k$, and then use Construction 1.5 to glue the second vertex of $G_c$ to the vertex $v_\ell$.

**Construction 1.12** *Dividing by $b$ each edge of a simple h-marked graph $G$.*

For each edge $e$ linking two vertices $v$ and $v'$ in the given graph $G$, perform the following operations. Remove $e$. Then introduce $b - 1$ new vertices $v_1(e), \ldots, v_{b-1}(e)$. Link with one edge: $v$ to $v_1(e)$, then $v_i(e)$ to $v_{i+1}(e)$ for $i = 1, \ldots, b-2$, and $v_{b-1}(e)$ to $v'$. Denote by $G'$ the resulting graph on $n + m(b-1)$ vertices. (Recall that $m$ is the number of edges of $G$.) Let $\{v_i, v_j\}$ denote the given pair of order $h$ for $G$, with associated vector $^tS = (s_1, \ldots, s_n)$. The pair $\{v_i, v_j\}$ also forms a pair of order $h$ in $G'$ when $G'$ is endowed with the following marking. Set the weight of $v_k$ to be $s_k b$. If $e$ is an edge between $v_k$ and $v_\ell$, set the weight of $v_r(e)$ to be $bs_k + (s_\ell - s_k)r$, for $r = 1, \ldots, b-1$. The reader will check that the resulting graph $G'$ is indeed $h$-marked. For a description of $\Phi(G')$ in terms of $\Phi(G)$, the reader may consult [Lor1], p. 281.

**Construction 1.13** *Coalescing an h-marked graph $G$ with an $h'$-marked graph $G'$*

Let $\ell := \mathrm{lcm}(h, h')$ and $g := \gcd(h, h')$. We can construct an $\ell$-marked graph $G''$ as follows. Let $v_1, \ldots, v_n$, denote the vertices of $G$, with $\{v_1, v_n\}$ the given pair of order $h$. Let $^tS = (s_1, \ldots, s_n = 0)$ be the associated vector, and assume that $s_i \geq 0$, for all $i$. Let $w_1, \ldots, w_{n'}$ denote the vertices of $G'$, with $\{w_1, w_{n'}\}$ the given pair of order $h'$. Let $^tS' = (0, s'_2, \ldots, s'_{n'})$ be the associated vector, and assume that $s'_i \leq 0$, for all $i$.

Let $G''$ denote the union of $G$ and $G'$, with the vertex $v_n$ glued to the vertex $w_1$. Number the vertices of $G''$ as $\{v_1, \ldots, v_n = w_1, w_2, \ldots w_{n'}\}$. Then $\{v_1, w'_{n'}\}$ is a pair of order $\ell$. Its associated vector is

$$^tS'' = (h's_1/g, \ldots, h's_{n-1}/g, 0, hs'_2/g, \ldots, hs'_{n'}/g).$$

Note that when $h = h'$, this construction glues together two $h$-marked graphs to produce a third one. In particular, gluing in such a way several graphs $G_h$ produces a graph where each pair of vertices has order $h$.

## 2  A description of all $h$-marked graphs

**Theorem 2.1** *Let $G$ be an h-marked graph. Let $\{v_i, v_j\}$ denote the given pair of order $h$, and let ${}^tS = (s_1, \ldots, s_n)$ denote the associated vector. Assume that $s_i < s_j$. Perform the following operations on $G$:*

1. *For each $k = 1, \ldots, n$, remove all edges between the vertices of weight $s_k$, as in Construction 1.9.*

2. *In the graph resulting from 1), glue together all vertices of the same weight, as in Construction 1.5.*

3. *In the graph resulting from 2), thicken all edges that can be thickened, as in Construction 1.6.*

4. *In the graph resulting from 3), glue together all vertices of the same weight, as in 1.5.*

*Then the graph resulting from 4) is the following graph $G'$:*

*where the vertex $v_i$ in $G$ corresponds to the vertex of weight $s_i$ in $G'$ and similarly, the vertex $v_j$ in $G$ corresponds to the vertex of weight $s_j$ in $G'$.*

*Conversely, any h-marked graph $G$ can be constructed using the following sequence of operations. Start with a string of graphs $G_h$ of length $\ell$ as in the graph $G'$ above (for some $\ell$ depending on $G$). Unglue various vertices of $G'$ to obtain a graph $G_1$. Perform inverse thickenings of some configurations of edges of $G_1$ to obtain a graph $G_2$. Unglue vertices of $G_2$ that were introduced by the inverse thickenings to obtain a graph $G_3$. Finally, add various graphs to $G_3$ to obtain $G$.*

*Proof:* Let us start with a graph $G$ with a pair $\{v_i, v_j\}$ of order $h$. Assume that the weight of $v_i$ is 0, and that the weight $s_j$ of $v_j$ is positive. Lemma 1.1 shows that all other weights on the given marking are non-negative and at most equal to $s_j$. Step 1) produces a graph in which two vertices of same weight are not linked by an edge. Thus it is possible to glue such vertices and perform Step 2). The resulting graph contains at most one vertex of any given weight. Step 3) requires performing a thickening of any edges linking two vertices whose weights differ by at least 2. After Step 3) is performed, the difference between the weights of any two adjacent vertices is 1. Call $\overline{G}$ the graph obtained after Step 3).

Consider the vertex $v_i$ in $\overline{G}$. Let us call $w_1, \ldots, w_r$ the vertices adjacent to $v_i$ in $\overline{G}$. Let $c_k$ denote the number of edges between $v_i$ and $w_k$. By definition of a marking, the pair $\{v_i, v_j\}$ has order $h = \sum_{k=1}^{r} c_k$. Step 4) requires us to glue together all vertices of a given weight in $\overline{G}$. We let the reader check that the only vertices of weight 1 in $\overline{G}$ are the vertices $w_1, \ldots, w_r$. Thus, after gluing together these $r$ vertices of weight 1, the resulting graph has a single vertex $v_i$ of weight 0 and a single vertex $x_1$ of weight 1. These two vertices are linked by $h$ edges. We may repeat the above argument starting with the vertex $x_1$ instead of $v_i$. After gluing the vertices of weight 2 adjacent to $x_1$, we obtain a graph with a single vertex of weight 2, linked to $x_1$ by $h$ edges. This process ends when all vertices of weights $s_j - 1$ are glued together. The resulting graph is a string of graphs $G_h$ of length $s_j$. We leave it to the reader to check that the last statement in Theorem 2.1 holds.

We turn now our attention to the question of giving explicit graph theoretical descriptions of pairs of vertices of order $h$. Recall that a *path* $\mathcal{P}$ on a graph $G$ is an ordered sequence of vertices and edges $v_0, e_1, v_1, e_2, \ldots, e_k, v_k$, such that the edge $e_i$ is linked to $v_{i-1}$ and $v_i$. We do not require that the vertices of the path $\mathcal{P}$ be distinct. On the other hand, we shall always assume that the edges of a path are all distinct. We shall say that the above path $\mathcal{P}$ has length $k$.

Let $\mathcal{P}_1$ and $\mathcal{P}_2$ be two paths from a vertex $v_i$ to a vertex $v_j$ in a graph $G$. Let us say that $\mathcal{P}_1$ and $\mathcal{P}_2$ *do not overlap* if they do not contain a common edge of $G$.

Let $\mathcal{P}_1, \ldots, \mathcal{P}_h$, be a set of $h$ paths of same length $\ell$ between $v_i$ and $v_j$. Number the edges of each path consecutively, so that the $h$ edges linked to $v_i$ receive the numbering 1 and the $h$ edges linked to $v_j$ are numbered $\ell$. We shall say that these paths form a *system of $h$ paths between $v_i$ and $v_j$* if, for each $1 \leq k \leq \ell$, the graph obtained from $G$ by removing all edges of $\mathcal{P}_1, \ldots, \mathcal{P}_h$ numbered $k$ is a disconnected graph. Note that this disconnected graph has exactly two components, the component of $v_i$ and the component of $v_j$. In particular, the endpoints of an edge numbered $k$ are not in the same component. Thus, a path $\mathcal{P}_r$ in a system of paths cannot pass twice through the same vertex.

Let us number consecutively the vertices on a path $\mathcal{P}_r$ of length $\ell$, with $v_i$

receiving the number 0 and $v_j$ receiving the number $\ell$. We claim that if two paths $\mathcal{P}_r$ and $\mathcal{P}_t$ intersect in a vertex $v$, and if $v$ is numbered $\alpha$ on the first path and $\beta$ on the second path, then $\alpha = \beta$. Indeed, suppose that $\alpha < \beta$. The vertex $v$ is the endpoint on $\mathcal{P}_r$ of an edge $e_r$ also numbered $\alpha$, and $v$ is the endpoint on $\mathcal{P}_t$ of an edge $e_t$ also numbered $\beta$. When removing all edges numbered $\beta$, the endpoints of the edge $e_t$ are in different connected components. This is not possible if $\alpha < \beta$, because in this case one can go from one endpoint of the edge $e_t$ to the other by first going back to $v_i$ along $\mathcal{P}_t$, and then following $\mathcal{P}_r$. Hence, $\alpha = \beta$. This argument also shows that if two paths pass through the same edge (and, hence, pass through both endpoints of the edge) then this edge has the same numbering for both paths.

We claim that each path $\mathcal{P}_i$ in a system of paths is a shortest path between $v_i$ and $v_j$. Indeed, suppose that the distance $d$ between $v_i$ and $v_j$ is smaller than the length $\ell$ of $\mathcal{P}_i$. Let $\mathcal{P}$ be a path of length $d$ between $v_i$ and $v_j$. Remove from $G$ all the edges numbered $k$ in the system of paths. Since the resulting graph is disconnected, at least one of these edges belongs to $\mathcal{P}$. Since the numbering of an edge is independent of the path, it is possible, if $d < \ell$, to find an integer $k$ with $1 \leq k \leq \ell$ such that none of the edges numbered $k$ in the given system of paths belong to $\mathcal{P}$. Then removing the edges numbered $k$ does not disconnect the graph $G$, a contradiction.

**Lemma 2.2** *Let $G$ be a connected graph with a system of $h$ non-overlapping paths between two vertices $v_i$ and $v_j$. Then the pair $\{v_i, v_j\}$ is a pair of order $h$ in $G$.*

*Proof:* To prove the lemma, it is sufficient to describe the marking of $G$ that exhibits $\{v_i, v_j\}$ as a pair of order $h$. Set the weight of $v_i$ to be 0, and the weight of $v_j$ to be $\ell$. Number the vertices of each path $\mathcal{P}_k$ consecutively, with $v_i$ numbered 0 and $v_j$ numbered $\ell$. As we noted above, this is a well-defined way of numbering the vertices of the system of paths. Then set the weight of the vertex numbered $r$ to be $r$. (Note that these weights on the system of paths is a marking of the system of paths.) Remove from $G$ all edges belonging to the paths $\mathcal{P}_k$, $k = 1, \ldots, h$. The new graph is the union of connected components. Each connected component $\mathcal{C}$ contains a vertex of the system of paths. By definition of a system of paths, if a connected component $\mathcal{C}$ contains a vertex in a path of the system with weight $r$, then all other vertices of the system of paths that are contained in $\mathcal{C}$ have also weight $r$. Thus, set the weight of any vertex on $\mathcal{C}$ to be $r$. The reader will check that this system of weights on the vertices of $G$ is indeed a marking of $G$. Note that we may now say that $G$ is obtained from the system of paths by addition of graphs (as in 1.4).

Note that the converse of Lemma 2.2 does not hold if $h \geq 3$. For instance, a cycle on $h \geq 3$ vertices does not contain a system of $h$ paths, but the cycle can be $h$-marked. The following corollary to Theorem 2.1 was first proven by Hans Gerd Evertz.

**Corollary 2.3** *Let $G$ be a connected graph. A pair $\{v_i, v_j\}$ has order 1 if and only if there is a system of one path in $G$ linking $v_i$ to $v_j$. In other words, a pair $\{v_i, v_j\}$ has order 1 if and only if there is a path $\mathcal{P}$ between $v_i$ and $v_j$ with the property that for all edges $e$ of $\mathcal{P}$, $G \setminus \{e\}$ is disconnected.*

*Thus $G$ is multiply connected if and only if its laplacian is spread.*

*Proof:* Lemma 2.2 shows that given a pair of vertices linked by such a path, then this pair has order 1. To prove the converse, let $\{v_i, v_j\}$ be a pair of order 1. Theorem 2.1 shows that this pair can be obtained by starting with a string on $\ell$ vertices (for some $\ell$) and performing certain operations. The first operation is the ungluing of vertices. Since the string is a simple graph, there are no possibilities of ungluing vertices. Similarly, there are no possibilities of performing the second operation, inverse thickening. Thus, the only possible operation is to add graphs. Hence, the resulting graph contains a path between $v_i$ and $v_j$, namely, the initial string, and removing any vertex of this string disconnect the graph by construction.

**Corollary 2.4** *Let $G$ be a connected graph. A pair $\{v_i, v_j\}$ has order 2 if and only if there is a system of 2 paths between $v_i$ and $v_j$.*

*Proof:* Since we do not assume that the two given paths between $v_i$ and $v_j$ do not overlap, we need to slightly modify the proof of Lemma 2.2 to show that the given pair of vertices has order 2. Number as usual the vertices of the paths consecutively. Define the weights of the vertices of the two paths using the following procedure: set the weight of $v_i$ to be zero. Suppose that a vertex numbered $k-1$ has been given the weight $s$. This vertex is linked to an edge $e$ numbered $k$. If $e$ belongs only to one of the two paths of the system, set the weight of the other endpoint of $e$ to be $s+1$. If $e$ belongs to both paths, set the weight of the other vertex of $e$ to be $s + 2$. Now number the vertices of $G$ that are not on the two given paths as in the proof of Lemma 2.2. We leave it to the reader to check that these weights form a well-defined marking of $G$, and that the pair $\{v_i, v_j\}$ has order 2.

To prove the converse, let $\{v_i, v_j\}$ be a pair of order 2. Theorem 2.1 shows that this pair can be obtained by starting with a string of $\ell$ graphs $G_2$ (for some $\ell$) and performing certain operations: First, the ungluing of certain vertices, and then inverse thickenings of edges. The resulting graph is a system of two paths between $v_i$ and $v_j$, and these two paths overlap if and only if an inverse thickening has been performed. The last operation listed in Theorem 2.1 is the addition of graphs. Thus, in the graph $G$, $v_i$ and $v_j$ are linked by a system of two paths.

**Example 2.5** A graph $G$ consisting in two vertices $v$ and $w$ linked by $h$ chains of length $\ell$ is a system of $h$ non-overlapping paths between $v$ and $w$. Thus $\{v, w\}$ is a pair of order $h$. The graph $G$ can also be obtained by dividing by $\ell$ each edge of the graph $G_h$.

Let us now discuss a slightly more general class of graphs which arise in the context of reduction of modular curves. Consider the graph $G$ consisting of two vertices $v$ and $w$, linked by $d$ chains of length $n_1, \ldots, n_d$, respectively.

Counting the number of spanning trees of $G$ is easy and we find that the order of the group $\Phi(G)$ is equal to $(\prod_{i=1}^{d} n_i)(1/n_1 + \ldots + 1/n_d)$. The structure of the group $\Phi(G)$ is discussed in [BLR], p. 283, and in [Ray], p. 18. (Note that in [Ray], the formula has typos. The correct formula is the one in [BLR].)

**Claim 2.6** *The pair $\{v, w\}$ has order $\mathrm{lcm}(n_1, \ldots, n_d)(1/n_1 + \ldots + 1/n_d)$.*

*Proof:* It is sufficient to exhibit the corresponding marking on the graph $G$. Let $L = \mathrm{lcm}(n_1, \ldots, n_d)$.

**Claim 2.7** *The pair* $\{v, v_{1,d}\}$ *has order* $\mathrm{lcm}(n_1, \ldots, n_{d-1})n_d(\sum_{i=1}^{d} 1/n_i)$.

More generally, the pair $\{v, v_{k,d}\}$ has order described as follows. Let $P$ be the smallest integer such that $P/n_1, \ldots, P/n_{d-1}$ and

$$Q := P + \frac{(n_d - k)P(1/n_1 + \ldots + 1/n_{d-1})}{k}$$

are all integers. In particular, $\mathrm{lcm}(n_1, \ldots, n_{d-1})$ divides $P$.

**Claim 2.8** *The order of the pair* $\{v, v_{k,d}\}$ *is* $(P/k)n_d(\sum_{i=1}^{d} 1/n_i)$.

Again, to prove these claims, we exhibit the appropriate marking:

The order of the pair $\{v, v_{k,d}\}$ is then

$$\begin{aligned}
h &= P/n_1 + \ldots + P/n_{d-1} + Q \\
&= (\sum_{i=1}^{d-1} 1/n_i)P + P(k + n_d - k)/kn_d + \frac{n_d - k}{k}P(\sum_{i=1}^{d-1} 1/n_i) \\
&= (\sum_{i=1}^{d} 1/n_i)Pn_d/k.
\end{aligned}$$

**Remark 2.9** The matrix $M$ endows the group $\Phi(G)$ with a canonical non-degenerate pairing $\langle\,,\,\rangle$ with values in $\mathbb{Q}/\mathbb{Z}$ (see [Lor5]). In particular, given an element $\tau$ of $\Phi(G)$, image of $E_{ij}$ modulo $\text{Im}(M)$, let $S$ be the associated vector with $MS = hE_{ij}$. Then $\langle \tau, \tau \rangle = (s_i - s_j)/h$. Thus the explicit computations of vectors $S$ in the above example can be used, for instance, to compute that for the element $\tau$ associated to the pair $\{v, w\}$ in 2.6, $\langle \tau, \tau \rangle = -(1/n_1 + \ldots + 1/n_d)^{-1}$.

**Remark 2.10** While the smallest $h$-marked graph is easy to describe (i.e., it is the graph $G_h$), the smallest simple graph(s) that can be $h$-marked are not know in general. It would also be interesting to know what is the largest prime number that can divide the number of spanning trees of a simple graph on $n$ vertices.

Clearly, starting with $G_h$ and dividing each edge in 2 produces a simple $h$-marked graph with $h+2$ vertices. The cycle on $h$ vertices can be $h$-marked. But, in general, it is possible to find simple $h$-marked graphs on $n$ vertices with $n$ much smaller than $h$. For instance, an odd integer $h \geq 9$ can be written in the form $h = xy + x + y$ with $x, y \geq 2$ if and only if $h$ is not of the form $h = 2p - 1$ with $p$ prime. For $h = xy + x + y$, we can use Example 2.6 with $n_1 = 1$, $n_2 = x$, and $n_3 = y$. This graph has $x + y$ vertices, and a pair of order $x + y + xy$. Thus in this example, $n \leq h/2$. Note that this method does not always produce, for a given $h$, graphs with the smallest ratio $n/h$. Indeed, with $h = 137 = 5 + 22 + 110$, we find that $n = 27$. The graph with $n_i = i$, $i = 1, \ldots, 5$, has a pair of order $h = 137$ and only $n = 12$ vertices.

When $h = 2p - 1$ with $p$ prime, we can use Example 2.6 with $n_1 = 1$, $n_2 = 3$, and $n_3 = \frac{1}{2}(3p - 1) - 1$. This graph has a pair of order $h$ and $n = \frac{3}{4}h + \frac{9}{4}$ vertices.

Note that if a pair $\{v_i, v_j\}$ is of order $h$ which is maximal among all possible orders of pair of vertices in simple graphs of size $n$, then $v_i$ and $v_j$ are linked by an edge. Indeed, if they are not linked by an edge, then consider the graph $G'$ obtained by adding an edge between $v_i$ and $v_j$. Suppose that in the original marking of $G$, $s_i = 0$ and $s_j = s > 0$. Then the pair $\{v_i, v_j\}$ has order $h + s > h$ in $G'$.

**Remark 2.11** It is possible to find examples of graphs $G$ that cannot be $h$-marked, even though the group $\Phi(G)$ contains an element $\varphi$ of order $h$. For instance, consider

It is easily checked that $|\Phi(G)| = 12$. Since $\Phi(G)$ has an element of order 12, we find that $\Phi(G) \cong \mathbb{Z}/12\mathbb{Z}$. The reader will check that $G$ has no pairs of order 2 or 3.

When $\Phi(G)$ is not cyclic, let $\varepsilon(G)$ denote the exponent of $\Phi(G)$. Since the vectors $E_{ij}$ form a basis of the orthogonal space to the vector $(1, 1, \ldots, 1)$, we find that the set $\{E_{ij}, i \neq j\}$ is a set of generators for $\Phi(G)$. It would be interesting to be able to determine directly from the graph whether some elements $E_{ij}$ has maximal order $\varepsilon(G)$. In all examples we computed, we always found a pair of vertices of $G$ of maximal order. Since the size of the graphs for which we were able to make the computation was relatively small, it is quite possible that an example of a graph $G$ where no pair of vertices has order $\varepsilon(G)$ could be found.

## 3 Collapsed values of matrices

Let us call a matrix $M \in M_n(\mathbb{Z})$ $\mu$-*collapsed*, or *collapsed for an (integer) value* $\mu$, if the matrix $M - \mu \text{Id}$ is not spread. In other words, $M$ is $\mu$-collapsed if the quotient map $\mathbb{Z}^n \to \mathbb{Z}^n/\text{Im}(M - \mu \text{Id})$ is not injective when restricted to the standard basis vectors $\{e_1, \ldots, e_n\}$. We shall call an integer $\mu$ such that $M$ is $\mu$-collapsed a *collapsed value of* $M$. Since $\det(M - x\text{Id})$ is the characteristic polynomial of $M$, any integer $\mu$ close enough to an eigenvalue of $M$ so that $|\det(M - \mu \text{Id})| < n$ is automatically a collapsed value of $M$. We shall see below that $\mu$ can be a collapsed value for $M$ even when $\det(M - \mu \text{Id}) \geq n$. Moreover, contrary to the case of the eigenvalues of $M$, the number of collapsed values of $M$ is not bounded by a constant depending on $n$ only.

Before considering the general case of any matrix $M$, let us first use Theorem 2.1 to study the collapsed values of a laplacian.

**Theorem 3.1** *Let $G$ be a connected simple graph. If $\mu \notin [0, n+1]$, then the laplacian of $G$ is not $\mu$-collapsed.*

*Proof:* Assume first that $\mu \geq n+2$. Consider $G$ as a subgraph of $K_\mu$, the complete graph on $\mu$ vertices. Let $C$ denote the complement of $G$ in $K_\mu$. In other words, $C$ is the subgraph obtained from $K_\mu$ by removing all edges of $G$. The matrix $M(C)$ has the following form. Let $J_\ell$ denote the $\ell \times \ell$ matrix all of whose entries are $+1$. Then

$$M(C) = \begin{pmatrix} & & & 1 & \cdots & 1 \\ & J_n - \mu\mathrm{Id} - M(G) & & \vdots & & \vdots \\ & & & 1 & \cdots & 1 \\ 1 & \cdots & 1 & -(\mu-1) & \cdots & 1 \\ \vdots & & \vdots & \vdots & \ddots & \vdots \\ 1 & \cdots & 1 & 1 & \cdots & -(\mu-1) \end{pmatrix}$$

Add all rows to the last row. Add all columns to the last column. Substract the penultimate row to all other rows. Then add all columns to the penultimate column. Finally, add the penultimate column to all but the last column to obtain a matrix $M'$ as follows:

$$M' = \begin{pmatrix} & & & 0 & \cdots & 0 & 0 & 0 \\ & -\mu\mathrm{Id} - M(G) & & \vdots & & \vdots & \vdots & \vdots \\ & & & 0 & \cdots & 0 & 0 & 0 \\ 0 & \cdots & 0 & -\mu & \cdots & 0 & 0 & 0 \\ \vdots & & & \vdots & \ddots & \vdots & \vdots & \vdots \\ 0 & \cdots & 0 & 0 & \cdots & -\mu & 0 & 0 \\ 0 & \cdots & 0 & 0 & \cdots & 0 & -1 & 0 \\ 0 & \cdots & 0 & 0 & \cdots & 0 & 0 & 0 \end{pmatrix}$$

The reader will check that $\mathrm{Im}(M')$ contains $E_{ij}$ with $1 \leq i < j \leq n$ if and only if $\mathrm{Im}(M(C))$ contains $E_{ij}$. By construction, the graph $C$ is multiply connected. Hence, $M(C)$ is 0-spread by 2.3. Thus, $\mathrm{Im}(M')$ does not contain $E_{ij}$, for all $1 \leq i\ j \leq n$. It follows that $-\mu\mathrm{Id} - M(G)$ is 0-spread.

Let us now treat the case where $\mu \leq -1$. Our proof below does not require the hypothesis that $G$ is simple. To ease the notation, let $\nu := -\mu$, and consider the graph $\overline{G}$ given by

$$M(\overline{G}) = \begin{pmatrix} & & & \nu \\ & M(G) - \nu\mathrm{Id} & & \vdots \\ & & & \nu \\ \nu & \cdots & \nu & -n\nu \end{pmatrix}$$

Add all rows to the last row and all columns to the last column to obtain a new matrix $M'$ with null last row and last column. As above, we let the reader verify that $\text{Im}(M(\overline{G}))$ contains $E_{ij}$ with $1 \leq i < j \leq n$ if and only if $\text{Im}(M')$ contains $E_{ij}$. Since $\overline{G}$ is multiply connected, we find that $M(\overline{G})$ is 0-spread. Thus $M'$, and then $M(G) - \nu\text{Id}$, are also 0-spread.

**Example 3.2** The reader will verify that the complete graph $K_n$ is $\mu$-collapsed if and only if $\mu = n - 1$ or $\mu = n + 1$. Thus the upperbound "$\mu \geq n + 2$ implies that $G$ is not $\mu$-collapsed" is sharp in Theorem 3.1. Note that an easy computation shows that

$$\mathbb{Z}^n/\text{Im}(M(K_n) - \mu\text{Id}) \cong (\mathbb{Z}/(\mu+n)\mathbb{Z})^{n-2} \times \mathbb{Z}/(\mu+n)\mu\mathbb{Z}.$$

The reader may have noted that the interval $[0, n]$ is the smallest interval that contains all the eigenvalues of all simple graphs $G$ of size $n$ (see, e.g., [Mer], p. 148). As we shall now discuss, a statement similar to Theorem 3.1 is true for a class of matrices larger than the class of the laplacians of simple graphs. Our next proposition generalizes the statement of Theorem 3.1 to any matrix $M \in M_n(\mathbb{Z})$. Proposition 3.4 applies to quadratic forms and provides an alternate proof of Theorem 3.1. Recall that if $v$ is any vector in $\mathbb{R}^n$, then $\|v\|$ denotes its euclidian length. Recall also that if $M$ is any matrix, $\|M\|$ denotes the maximal value taken by the function $\|Mv\|$ when $v$ is on the unit sphere.

**Proposition 3.3** Let $M \in M_n(\mathbb{Z})$. Assume that $M$ is $\mu$-collapsed. Then $|\mu| \leq \|M\| + \sqrt{2}$. If the column vectors of $M - \mu\text{Id}$ do not contain a vector $E_{ij}$, then $|\mu| \leq \|M\| + 1$.

*Proof:* Without loss of generality, we may assume that there exists a vector $v$ with transpose ${}^t v = (v_1, \ldots, v_n)$ such that $(M - \mu\text{Id})v = E_{12}$. Thus

$$\|E_{12} + \mu v\| \leq \|M\|\,\|v\|.$$

Hence,
$$2 + 2\mu(v_1 - v_2) \leq (\|M\|^2 - \mu^2)\|v\|^2.$$

Using Lagrange multipliers, we find that $v_1 - v_2 \geq -\sqrt{2}\|v\|$. Thus

$$(\|M\|^2 - \mu^2)\|v\|^2 \geq 2 + 2\mu(v_1 - v_2) \geq 2 - 2|\mu|\sqrt{2}\|v\|.$$

It follows that $|\mu| \leq \sqrt{2}/\|v\| + \|M\|$. If the column vectors of $M - \mu\text{Id}$ do not contain $E_{ij}$, for all $i \neq j$, we find that $v \neq e_i$, for all $i$. Hence, $\|v\| \geq \sqrt{2}$, and $|\mu| \leq 1 + \|M\|$.

**Proposition 3.4** *Let $M$ be a semidefinite positive quadratic form. Denote by $\lambda_1$ its largest eigenvalue.*

*(i) If $\mu \notin [-1, \lambda_1 + \sqrt{2}]$, then $\mu$ is not a collapsed value of $M$.*

*(ii) The value $\mu = -1$ is a collapsed value of $M$ if either two columns of $M$ are equal or if $Me_i = -e_j$ for some $i \neq j$.*

*(iii) If $M$ is positive definite or is the laplacian of a connected graph, then $\mu < 0$ is not a collapsed value of $M$.*

*(iv) When $M$ is positive definite, let $\lambda_n > 0$ denote its smallest eigenvalue. If $\mu$ is a collapsed value, then $\mu \geq \lambda_n - \sqrt{2}$.*

*Proof:* When $M$ is a semidefinite positive quadratic form, $\|M\| = \lambda_1$. Thus, we may apply Proposition 3.3 to show that $\mu > \lambda_1 + \sqrt{2}$ is not a collapsed value. Now assume that $\mu < 0$ and let ${}^t v = (v_1, \ldots, v_n)$ be an integer vector such that $(M - \mu\text{Id})v = E_{ij}$. Then

$$ {}^t v(M - \mu\text{Id})v = {}^t vMv - \mu\|v\|^2 = v_i - v_j, $$

and thus ${}^t vMv + (|\mu|-1)\|v\|^2 = (v_i - v_j) - \|v\|^2$. The left side of this equality is non-negative by hypothesis. It is easy to check that the right side is non-positive. Thus for this equality to hold, both sides have to be 0. Hence, $\mu = -1$ and part (i) is proved.

In addition to $\mu = -1$, for both sides to be 0, the reader will check that $v$ must either equal $e_i$, or $-e_j$, or $E_{ij}$. When such is the case, we find that $Me_i = -e_j$, or $ME_{ij} = 0$. This proves part (ii). To prove part (iii), note that if two columns of $M$ are equal, then $M$ cannot be definite positive and cannot be a laplacian. Moreover, if $Me_i = -e_j$, then ${}^t e_i Me_i = 0$, and $M$ cannot positive definite. If $M$ is a laplacian, then $M$ cannot have a column with only a single non-zero entry. Thus (iii) holds. Part (iv) is similar to the proof of 3.3.

Note that Theorem 3.1 follows easily from 3.4 and from the fact that the largest eigenvalue of the laplacian of a simple graph is at most $n$. Indeed, a collapsed value $\mu$ is an integer at most equal to $\lambda_1 + \sqrt{2}$. Thus $\mu \leq n + 1$.

There are other useful upper bounds for the largest eigenvalue of the laplacian of a simple graph. For instance, if $d_1 \geq d_2 \geq \ldots \geq d_n$ is the ordered sequence of the degrees of vertices of $G$, then it is known that $\lambda_1 \leq d_1 + d_2$ (see, e.g., [Mer], p. 148). Consider for instance a string of $n \geq 4$ vertices. By hypothesis, $d_1 + d_2 = 4$. The reader may check that the laplacian of such a string is $\mu$-collapsed only for $\mu = 0, 1, 2$, and $3$. In the case of a bipartite graph $K_{p,q}$, with $p, q \geq 2$, the reader may check that the only collapsed values of the laplacian of $K_{p,q}$ are $\mu = p \pm 1$ and $\mu = q \pm 1$.

**3.5** We conclude this section with examples of matrices having a large number of collapsed values. Note first that if $E_{ij}$ is an eigenvector for an eigenvalue $\lambda$ of $M$, then $\lambda$ is an integer. Moreover, $M$ is then $\mu$-collapsed if $\mu = \lambda \pm 1$. Indeed, $ME_{ij} = \lambda E_{ij}$ implies that $(M - (\lambda - 1)\text{Id})E_{ij} = E_{ij}$ and $(M - (\lambda + 1)\text{Id})(-E_{ij}) = E_{ij}$.

Let $n \geq 8$ be even. We can use this remark to produce examples of matrices $M(G)$ where the number of collapsed values $\mu$ is equal to $1 + n/2$. Indeed, consider the graph on $n$ vertices given by

$$M(G) = \begin{pmatrix} -(n-1) & 1 & 1 & 1 & 1 & 1 & \cdots & 1 & 1 & 1 & 1 \\ 1 & -(n-1) & 1 & 1 & 1 & 1 & \cdots & 1 & 1 & 1 & 1 \\ 1 & 1 & -(n-4) & 0 & 1 & 1 & \cdots & 1 & 1 & 0 & 0 \\ 1 & 1 & 0 & -(n-4) & 1 & 1 & \cdots & 1 & 1 & 0 & 0 \\ 1 & 1 & 1 & 1 & * & 0 & \cdots & 0 & 0 & 0 & 0 \\ 1 & 1 & 1 & 1 & 0 & * & & 0 & 0 & 0 & 0 \\ \vdots & \vdots & \vdots & \vdots & \vdots & & \ddots & \vdots & \vdots & \vdots \\ 1 & 1 & 1 & 1 & 0 & 0 & & -4 & 0 & 0 & 0 \\ 1 & 1 & 1 & 1 & 0 & 0 & \cdots & 0 & -4 & 0 & 0 \\ 1 & 1 & 0 & 0 & 0 & 0 & \cdots & 0 & 0 & -2 & 0 \\ 1 & 1 & 0 & 0 & 0 & 0 & \cdots & 0 & 0 & 0 & -2 \end{pmatrix}$$

It is easy to check that the rank of $M(G) - \lambda \text{Id}$ is smaller than $n$ when $\lambda = -n, -(n-4), -(n-6), \ldots, -4, -2, 0$. Thus, using the above remark, we find that $M(G)$ is $\mu$-collapsed for $\mu = -(n+1), -(n-1), -(n-3), -(n-5), \ldots, -3, -1$. Since $G$ is multiply connected, $M(G)$ is not 0-collapsed. It would be interesting to determine the maximal number of collapsed values that a graph on $n$ vertices can have.

With matrices $M$ that are not laplacians of graphs, one can find examples with considerably more collapsed values $\mu$. Consider for instance

$$M = \begin{pmatrix} -a & 1 \\ 1 & a \end{pmatrix}.$$

The eigenvalues of $M$ are $\lambda = \pm\sqrt{a^2 + 1}$. The reader will check that $M$ is $\mu$-collapsed for $\mu = a, a-1, -a, -(a+1)$. For specific $a$, it is sometimes possible to lengthen this list. For instance, when $a = 1$, we find that $M$ is $\mu$-collapsed for $\mu = 2, 1, 0, -1, -2$. Let $f(x) = x^2 - 2$ be the characteristic polynomial of $M$. Note that even though $f(2) = |f(0)| = f(-2) \geq n$, the values $2, 0,$ and $-2$ are collapsed values of $M$.

Consider the matrix $M_k$ built as a diagonal block matrix, with $k$ $(2 \times 2)$-blocks given as follows:
$$M + 5\text{Id}, \ldots, M + 5k\text{Id},$$
(with $a = 1$ in $M$). The eigenvalues of $M_k$ are $5i \pm \sqrt{2}$, $i = 1, \ldots, k$. It is easy to check that $M_k$ is $\mu$-collapsed for all integer values $\mu$ with $3 \leq \mu \leq 5k+2$. This example shows that Proposition 3.3 is sharp : $\|M_k\| = \lambda_1 = 5k + \sqrt{2}$ since $M_k$ is positive definite, and all integers $\mu$ such that $\lambda_n - \sqrt{2} \leq \mu \leq \|M_k\| + \sqrt{2}$ are collapsed values for $M_k$.

**Example 3.6** The following example, communicated to us by A. Granville and C. Pomerance, shows that the number of collapsed values of a matrix $M$ is not bounded by a constant depending on $n$ only. Consider
$$M = \begin{pmatrix} 2\ell & 1 \\ -\ell^2 & 0 \end{pmatrix}.$$
The characteristic polynomial of $M$ equals $(x - \ell)^2$. The solutions of the equation $(s_1, s_2)(^tM - \mu\text{Id}) = (1, -1)$ are given by $s_1 = (1-\mu)/(\mu-\ell)^2$, and $s_2 = (\ell^2 - 2\ell + \mu)/(\mu - \ell)^2$.

Fix any integer $r$ and a finite list of primes $p_1, \ldots, p_r$. The Chinese Remainder Theorem shows that it is possible to find an integer $\ell$ such that, for each $i = 1, \ldots, r$, $\ell \equiv p_i + 1 \pmod{p_i^2}$. It is easy to check that when $\mu = \ell - p_i$, the numbers $s_1$ and $s_2$ are integers. Thus, $M$ has at least $r$ collapsed values.

Note that the characteristic polynomial of $M$ is very special. It would be interesting to know whether a similar example exists with a more "generic" characteristic polynomial, having for instance two distinct roots which are not integers.

# References

[BLR]   S. Bosch, W. Lüktebohmert, and M. Raynaud, *Néron Models*, Springer Verlag, 1990.

[GMW]   R. Grone, R. Merris, and W. Watkins, *Laplacian Unimodular Equivalence of Graphs*, in Combinatorial and Graph-Theoretic Problems in Linear Algebra, R. Brualdi et al Eds., Springer Verlag, 1995.

[Lor1]   D. Lorenzini, *A finite group attached to the laplacian of a graph*, Discr. Math. **91** (1991), 277-282.

[Lor2]   D. Lorenzini, *Research Problem*, Discr. Math. **84** (1990), 327-329.

[Lor3]   D. Lorenzini, *Arithmetical graphs*, Math. Ann. **285** (1989), 481-501.

[Lor4]   D. Lorenzini, *Reduction of points in the group of components of the Néron model of a jacobian*, in preparation.

[Lor5]   D. Lorenzini, *Pairings on groups of components*, in preparation.

[Mer]   R. Merris, *Laplacian matrices of graphs: A Survey*, Linear Algebra and its Applications **197-198** (1994), 143-176.

[Ray]   M. Raynaud, *Jacobienne des courbes modulaires et opérateurs de Hecke*, in Astérisque **196-197** (1991), 9-25.